# On the Riemann zeta-function, Parts IV-V

**Part IV: On the Riemann zeta-function and meromorphic characteristic functions.**
**Part V: A relation of its nonreal zeros and first derivatives thereat to its values on ½ + 4N.**

## By Anthony Csizmazia

### E-mail: apcsi2000@yahoo.com

## Abstract

In Part I an odd meromorphic function f(s) has been constructed from the Riemann zeta-function evaluated at one-half plus s. The conjunction of the Riemann hypothesis and hypotheses advanced by the author in Part I is assumed. In Part IV we derive the two-sided Laplace transform representation of f(s) on the open vertical strip V of all s with real part between zero and four. An additional hypothesis is used to prove that the Laplace density of f(s) on the strip V is positive. Let z(n) be the nth critical zero of the Riemann zeta-function of positive imaginary part in order of magnitude thereof. In Part V an expression is derived for z(1). A relation is obtained of the pair z(n) and the first derivative thereat of the zeta-function to the preceding such pairs.

**MSC (Mathematics Subject Classification).** 11Mxx Zeta and L-functions: analytic theory. 11M06 $\zeta(s)$ and L(s, $\chi$). 11M26 Nonreal zeros of $\zeta(s)$ and L(s, $\chi$); Riemann and other hypotheses. 30xx Functions of a complex variable. 44A10 Laplace transform. 42A82 Positive definite functions. 60E10 Characteristic functions; other transforms.

## Table of contents.

**Abstract. Keywords. MSC (Mathematics Subject Classification).**

## Introduction

## Part IV: On the Riemann zeta-function and meromorphic characteristic functions.

## §1 The conditional Laplace representation of f(s) on $V_0$.

## §2 Proof of the Main conditional theorem (1).
Main conditional theorem (1)

(i) The equality of the conditional and unconditional Laplace densities.
(i') The boundedness of the density.
(ii) The conditional extension of the unconditional Laplace representation of f(s) on $V_0'$ to $V_0$.

## §3 Proofs of the Main conditional theorems (2)-(3).

## §4 Metric norms and analytic characteristic functions.

## References

## Part V: A relation of its nonreal zeros and first derivatives thereat to its values on ½ + 4N.

## Index of symbols

## §1 Proof of the conditional relations of $\gamma_n$, $\zeta'(½ + i\gamma_n)$ to their predecessors and the $\zeta(½ + 4k)$.

## §2 Representation of $p_{i,+}(z)$ via j(y).

## References

On the Riemann zeta-function, Part IV.
On the Riemann zeta-function and meromorphic characteristic functions.

(Complex plane) C. (Real line) **R**.

## §1 The conditional Laplace representation of f(s) on $V_0$.

**Review** Part I, §4, (4.3); §5, Introduction, (5.3) and (5.4).

**Conditional theorem 1.1** *Assume C^. On $V_0$: f(s) = $\int_R d(y)e^{sy}g_0(y)$.*

**Proof** The Conditional theorem 2.2 proven in Part III, §2, is stated as Conditional theorem 5.1 in Part I, §5, Introduction.(see A. Csizmazia [3, 5]). It established that C^ implies that f(s) is represented on C − Z° by its formal partial fraction expansion p(s). Thus f(s) = p(s) on $V_0$. C^ implies A is finite. Then the Conditional corollary 5.1 of Part I, §5, (see A. Csizmazia [3]), gives the asserted Laplace representation.

## §2 Proof of the Main conditional theorem (1).

**Review** Part I, §5, (5.4).

**Main conditional theorem (1)** *Assume C^.*
**(i) The equality of the conditional and unconditional Laplace densities.**
*If y is real, then:*

**Eq. (*)**

$$\lambda(y) + c(0) - P_0(\pi e^{2y}) = P_0(\pi e^{-2y}).$$

**(i′) The boundedness of the density.**
*$P_0(v)$ is bounded on the real axis.*
**(ii) The conditional extension of the unconditional Laplace representation of f(s) on $V_0′$ to $V_0$.**
*On the strip $V_0$:* $f(s) = \int_R d(y)e^{sy}P_0(\pi e^{-2y})$.

**Proof of (i).** $g_0(y) = P_0(\pi e^{-2y})$, for y > 0. Say s = x + it, with ½ < x < 4 and t real. Consider the previous Conditional theorem 1.1 (see also Part I, §5, (5.4)) and the Main unconditional theorem (stated in Part I, §3, and) proven in Part II, §6. (See A. Csizmazia [3-4].) Together they yield:
(″) $\int_{y<0} d(y)e^{ity}(e^{xy}\theta(y)) = 0$, with $\theta(y) := g_0(y) - P_0(\pi e^{-2y})$.

C^ implies A is finite.

**Conditional claim (1°)** *Assume A is finite. Fix x > ½. $e^{xy}|\theta(y)|$ vanishes with exponential rapidity as y recedes to -∞, with y < 0.*

Proof of Conditional claim (1°)

**Review** Theorem 3.2 and its consequence (′) of Part I, §3. The Unconditional theorem 6.1 (4) of Part II.

$e^{xy}|P_0(\pi e^{-2y})|$ vanishes with exponential rapidity as y recedes to -∞. $e^{xy}|g_0(y)|$ behaves likewise when x > 0, since ∣λ(y) ∣≤ A.

**Conditional claim (2°)** *Assume A < ∞. $e^{xy}\theta(y)$ is continuous in y.*

Proof of Conditional claim (2°) $P_0(z)$ is entire. Hence $P_0(\pi e^{-2y})$ is continuous. A is finite. So λ(y) is defined and continuous, for real y. Thus $g_0(y)$ is continuous. Hence so is $e^{xy}\theta(y)$.

The vanishing of the Fourier transform in (″) and the Conditional claims (1°), (2°) together imply $g_0(y) = P_0(\pi e^{-2y})$, for y < 0.

Proof of (ii). In the previous Conditional theorem 1.1 apply $g_0(y) = P_0(\pi e^{-2y})$ of (i) of this Main conditional theorem (1).

**Conditional continuity criterion** Assume $C^{\wedge}$.
$\lim y < 0,\ y \to 0\ g_0(y) = \lim y > 0,\ y \to 0\ g_0(y)$. Thus

$$\sum_{k \geq 1} c(i\gamma k) = -(c(0)/2 + \sum_{k \geq 1} c(4k)).$$

The previous criterion and the Conditional theorem 2.1 stated next follow from the Main conditional theorem (1), as observed in Part I, §5, (5.4). Therein Conditional theorem 2.1 is stated as Conditional theorem 5.3 (1).

**Conditional theorem 2.1**
*Assume $C^{\wedge}$. $g_0(z) = P_0(\pi e^{-2z})$, respectively*

$$\lambda(z) = -c(0) + P_0(\pi e^{2z}) + P_0(\pi e^{-2z}) = -(c(0) + 2\sum_{k \geq 1} c(4k)\cosh(4ky)),$$

*holds on the real line and so extends $g_0$, respectively $\lambda$, to an entire function on C of period $i\pi/2$.*

**§3 Proofs of the Main conditional theorems (2)-(3).**

Assume $C^{\wedge}$. Apply the Main conditional theorem (1), proven in §1. Also assume C5. Then (iii) Synthesis, Positivity of $P_0(v)$, Conditional Lemma 7.2 of §7, Part I, gives: $P_0(v) > 0$, for all $v > 0$; and $\inf_{v > \varepsilon} P_0(v) > 0$, when $\varepsilon > 0$. So the following Main conditional theorem (2), stated in Part I, §7, is attained.

**Main conditional theorem**
(2) *Assume $C^{\wedge}$ and C5. $f(s)$ is an analytic characteristic function on $V_0$:*

$$f(s) := 1/n(s) = \int_R d(y)e^{sy}P_0(\pi e^{-2y}),$$

*with $P_0(v)$ positive for $v > 0$. Also $\inf_{v > \varepsilon} P_0(v) > 0$, for any $\varepsilon > 0$.*

Apply the previous Main conditional theorem (2) together with the Main unconditional theorem (4), stated in Part I, §3, and proven in Part II, §6, The Mellin transform representation of $f(s, \beta)$, Results when $\beta = ¼$. As observed in Part I, §7, one obtains the Main conditional theorem (3) restated next.

**Main conditional theorem**

(3) *Assume C^ and C5. f(s) is a meromorphic characteristic function on the complex plane: When w is an integer and s is in $V_{4w}$,*

$$(-1)^w f(s) = \int_R d(y) e^{sy} P_{4w}(\pi e^{-2y}),$$

*with $P_{4w}(z)$ entire in z and $P_{4w}(v)$ positive for v > 0.*

Applying the relation f(-s) = -f(s) one obtains the counterparts of the above results for the negative half-plane of s with Re(s) < 0. When w ≤ -1, set $P_{4w}(\pi v)$ := $P_{-4(w+1)}(\pi/v)$, for v > 0.

## §4 Metric norms and analytic characteristic functions.

**Review** Part I, §3, A geometric consequence of the Main unconditional theorem (4), and §7, Metric norms and analytic characteristic functions. Part II, §6, Metric result when β = ¼ .

The next result emanates from the association of Corollary 2.2 of Part VI, Corollary 6.2 with β = ¼ of Part II and the Main conditional theorem (1) (ii) of §2 above.

**Conditional corollary 4.1** *Say Corollary 2.2 of Part VI holds. Assume C^ and C5. Let x, t be real with x not a multiple of four.*

*$m_x(t)$ is a metric norm in t on the real line.*

*$d_x(t_1, t_2)$ is a (finite-valued) translation invariant metric in $t_1$, $t_2$ on the real line.*

## On the Riemann zeta-function: Part V.
**A relation of its nonreal zeros and first derivatives thereat to its values on ½ + 4N.**

**Index of symbols**

$h^{\#}(z)$, L, D(r, ω), $D_1$(r, ω), S(z, ω), $S_{\leq}$(z, ω), Δ(z, ω), c(h, k, ω) and c(h, k).

**§1 Proof of the conditional relations of $\gamma_n$, $\zeta'$(½ + i$\gamma_n$) to their predecessors and the ζ(½ + 4k).**

**Review Part I §6**, **§3**, **Theorem 3.2 (i)**. **Main conditional theorem (1) (i)**, **§5**, **(5.4)**, **and Part IV**, **§2**.

In **Part I, §6**, the **Conditional corollary 6.1** is stated and its proof is deferred to this Part V. The **Conditional claim 6.1 (1)** was proven assuming (*) thereof and using **Conditional corollary 6.1 (1)**. It was observed that C^ and the **Main conditional theorem (1) (i)**, with y < 0, together imply (*). The **Conditional claim 6.1 (2)** resulted from the **Conditional claim 6.1 (1) and Conditional corollary 6.1 (2)**.

**Conditional corollary 6.1** of **Part I, §6,** is proven as **Conditional corollary 1** using **Lemma 1** below.

**Definition of the Poisson transform $h^{\#}(z)$ of h(y).** See (1), (2) of the next lemma.

**Lemma 1** *Assume ('): h(u) is an even entire function of u, $\theta < 1$ and $|h(y)| \leq K(y^\theta)$, for large positive y.*
*(1) $z \cdot \int_{y \geq 0} d(y)(1/(y^2 + z^2))h(y)$ converges absolutely to an analytic function, $h^\#(z)$, on the half-plane $Re(z) > 0$.*
*(2) $h^\#(z)$ extends to an entire function on C.*

**Proof** The proof of the **Lemma 1** is achieved by establishing a series of claims.

**Definitions of L, D(r, ω), D₁(r, ω).** Set L := $[-\infty, -1]U[1, \infty]$. Fix $\omega > 0$. Say $r > 0$. Let $D(r, \omega)$, $D_1(r, \omega)$ be the set of all u distant from S by at least r, with S respectively iωL, iω[-1, 1].

**Claim 1** *Say $\omega > 0$ and $y \geq \omega$. Assume (*): h(y) is a complex-valued measurable function, $\theta < 1$ and $|h(y)| \leq K(y^\theta)$, for $y \geq \omega$.*
*(1) $\int_{y \geq \omega} d(y)(1/(y^2 + z^2))h(y)$ converges uniformly absolutely in z on D(r, ω).*
*(2) $\int_{y \geq \omega} d(y)(1/(y^2 + z^2))h(y)$ is analytic in z on $C - i\omega L$.*

**Definition of S(z, ω).** Let S(z, ω) be the integral in (2) of Claim 1.

Proof of Claim 1. $1/|y^2 + z^2| = y^{-2} \cdot |1 + (z/y)^2|^{-1}$. Now $1/|1 + (z/y)^2|$ is bounded for z, y with z on D(r, ω) and $y \geq \omega$. Apply (*) of Claim 1 to obtain $\int_{y \geq \omega} d(y)|h(y)/(y^2 + z^2)| \leq K_1 \int_{y \geq \omega} d(y)y^{-2+\theta}$, for some $K_1$. Thus (1) and (2) hold.

Claim 2. *Taylor series for S(z, ω). Assume (*), of Claim 1, and $|z| < \omega$. S(z, ω) $= \sum_{k \geq 0} z^{2k}((-1)^k \cdot \omega^{-2k-1} \cdot \int_{y \geq 1} d(y)y^{-2(k+1)}h(\omega y))$, with each integral and the series converging absolutely.*

Proof of Claim 2. S(z, ω) $= \omega^{-1} \cdot \int_{y \geq 1} d(y)(1/(y^2 + u^2))h(\omega y)$, with u := z/ω. $|u| < 1$. $1/(y^2 + u^2) = \sum_{k \geq 0} (-1)^k \cdot u^{2k} \cdot y^{-2(k+1)}$. Consider $\sum_{k \geq 0} |u|^{2k} \cdot t(k, \omega)$, with t(k, ω) := $\int_{y \geq 1} d(y)y^{-2(k+1)}|h(\omega y)|$. Apply (*) to obtain t(k, ω) $\leq K \cdot \omega^\theta/(2k + 1 - \theta)$. So $\sum_{k \geq 0} |u|^{2k} \cdot t(k, \omega) \leq K \cdot \omega^\theta \cdot E(|u|)$, with $E(r) = \sum_{k \geq 0} r^{2k}/(2k + 1 - \theta)$ when $|r| < 1$. Say $-1 < r < 1$. Then $E(r) \leq 1/(1 - \theta) - \frac{1}{2}\log(1 - r^2)$.

Claim 3. *Say $\omega > 0$. Assume $\int_{0 \leq y \leq \omega} d(y)|h(y)|$ converges.*
*(1) $\int_{0 \leq y \leq \omega} d(y)h(y)/(y^2 + z^2)$ converges absolutely in z on D₁(r, ω).*
*(2) $\int_{0 \leq y \leq \omega} d(y)(1/(y^2 + z^2))h(y)$ is analytic in z on C - iω[-1, 1].*

**Definition of S≤(z, ω).** Let S≤(z, ω) be the integral in (2).

Proof of Claim 3 (1).
$1/(|z - iy| \cdot |z + iy|) \leq r^{-2}$ leads to $\int_{0 \leq y \leq \omega} d(y)|h(y)/(y^2 + z^2)| \leq r^{-2} \cdot \int_{y \leq \omega} d(y)|h(y)|$.

Say Re(u) > 0. Set q(u) := u·$\int_{0 \leq y \leq 1}$ d(y)/(y$^2$ + u$^2$). If y ≠ ±iz, then z/(y$^2$ + z$^2$) =(2i)$^{-1}$(1/(y − iz) − 1/(y + iz)). Thus ('): q(u) = (2i)$^{-1}$·$\sum_{\sigma = \pm 1}$ σ$\int_{[0, 1] - \sigma iu}$ d(y)/y.

Say u is in C − (−∞, 0].

Take log(u) to be the principal branch log(|u|) + iθ, with u= |u|·e$^{i\theta}$ and −π < θ < π. Take arctan(s) to be the principal branch that for real x has −π/2 < arctan(x) < π/2 and has iL as branch cut. arctan(s) = (2i)$^{-1}$·$\sum_{\sigma = \pm 1}$ σlog(1 + σiu). Let arccot(s) = π/2 − arctan(s).

Say Re(u) > 0. Evaluate q(u) using (') above. σ$\int_{[0, 1] - \sigma iu}$ d(y)/y = σlog(1 - σiu) + (iπ/2) − σlog(u). Note the singular behavior of log(u) as u approaches 0. Fortunately in q(u) the log(u) terms cancel one another, since the coefficient of their total contribution (2i)$^{-1}$·(−1)($\sum_{\sigma = \pm 1}$ σ) is 0. Thus

$$q(u) = \pi/2 + (2i)^{-1}\cdot\textstyle\sum_{\sigma = \pm 1} \sigma\log(1 - \sigma iu) = \operatorname{arccot}(u).$$

Hence q(u) has an analytic extension from the half-plane Re(u) > 0 to C − iL.

Claim 4. *Assume ω > 0 and h(u) is an even entire function of u. zS$_{\leq}$(z, ω) extends from the half-plane Re(z) > 0 to an analytic function on C − iωL.*

Proof of Claim 4.
Assume ω > 0 and Re(z) > 0. $\int_{0 \leq y \leq \omega}$ d(y)|(j(y$^2$) - j(-(z$^2$)))/(y$^2$ + z$^2$)| converges.

**Definition of Δ(z, ω).** Set Δ(z, ω) := $\int_{0 \leq y \leq \omega}$ d(y)(j(y$^2$) - j(-(z$^2$)))/(y$^2$ + z$^2$).

zS$_{\leq}$(z, ω) = z·Δ(z, ω) + h(iz)·arccot(z/ω).

h(s) = j(s$^2$), with j(u) an entire function of u. Set j$_k$ = j$^{(k)}$(0)/(k!).

Subclaim *Say ω ≥ 0.*
(1) *Assume that z is on C - iω[-1, 1].*

$$\Delta(z, \omega) = \textstyle\sum_{w \geq 0} (-(z^2))^w \cdot \sum_{k \geq w + 1} (j_k/(2(k-w)-1))\cdot \omega^{2(k-w)-1},$$

*with each of the series absolutely convergent.*
(2) *Δ(z, ω) is extended to an entire function of z, ω by the latter series.*

Proof of subclaim Say u ≠ v. (j(u) - j(v))/(u - v) = $\sum_{k \geq 1}$ j$_k$·(u$^k$ - v$^k$)/(u - v). Also (u$^k$ - v$^k$)/(u - v) = $\sum_{0 \leq w \leq k-1}$ u$^{k-1-w}$·v$^w$. Assume |u|, |v| ≤ B. Then $\sum_{0 \leq w \leq k-1}$ |u|$^{k-1-w}$·|v|$^w$ ≤ k·B$^{k-1}$. Let m(s) := $\sum_{k \geq 1}$ |j$_k$|·s$^k$. One has $\sum_{k \geq 1}$ |j$_k$|·k·B$^{k-1}$ = m'(B). Take v

= -(z²). Assume ω ≤ B^½ and |z| ≤ B^½. Then $\int_{y \le \omega} d(y) \sum_{k \ge 1} |j_k| \cdot \sum_{0 \le w \le k-1} |y^2|^{k-1-w} \cdot |v|^w \le m'(B) \cdot B^½ < \infty$. Thus the subclaim is valid.

$zS_{\le}(z, \omega) = z \cdot \Delta(z, \omega) + h(iz) \cdot arccot(z/\omega)$, for Re(z) > 0. That equality and Subclaim (2) yield Claim 4.

The Maclaurin expansion, for z with |z| < ω, of $zS_{\le}(z, \omega)$ is determined by the latter equality from the expansions of $\Delta(z, \omega)$ and $h(iz) \cdot arccot(z/\omega)$.

Corollary *Assume (') of **Lemma 1**. Let ω > 0. h^#(z) has an analytic continuation from the half-plane Re(z) > 0 to C − iωL.*

Proof of Corollary. If Re(z) > 0, then $h^\#(z) = zS(z, \omega) + zS_{\le}(z, \omega)$. That identity provides the analytic continuation of h^#(z) to C − iωL, since each of S(z, ω) and $zS_{\le}(z, \omega)$ is analytic there. So the corollary holds.

We now complete the proof of **Lemma 1**.

$$h^\#(z) = z(S(z, \omega) + \Delta(z, \omega)) + h(iz) \cdot arccot(z/\omega),$$

for z on C − iωL. So on C − iωL the odd part h₁(z) of h^#(z) is given by

$$h_1(z) := \tfrac{1}{2}(h^\#(z) - h^\#(-z)) = z(S(z, \omega) + \Delta(z, \omega)) - h(iz)arctan(z/\omega).$$

Apply Claim 2 together with the Subclaim to obtain the following. Say ω > 0. Each of $\int_{y \ge \omega} d(y)y^{-2(k+1)}h(y)$ and $\sum_{n \ge -k} (j_{n+k}/(2n-1)) \cdot \omega^{2n-1}$ converges absolutely.

**Definition of c(h, k, ω).**
Set $c(h, k, \omega) := (\int_{y \ge \omega} d(y)y^{-2(k+1)}h(y)) + \sum_{n \ge -k} (j_{n+k}/(2n-1)) \cdot \omega^{2n-1}$.

Say |z| < ω. $\sum_{k \ge 0} (z^{2k+1})(-1)^k \cdot c(h, k, \omega)$ converges absolutely. h₁(z) is analytic for z with |z| < ω: $h_1(z) = \sum_{k \ge 0} (z^{2k+1})(-1)^k \cdot c(h, k, \omega)$.

$c(h, k, \omega) = (-1)^k \cdot (h_0(z))^{(2k+1)}(0)/(2k+1)!$. Therefore c(h, k, ω) is constant in ω: $\partial_\omega(c(h, k, \omega)) = 0$.

**Definition of c(h, k).** Let c(h, k) be the common value of the c(h, k, ω) with ω positive.

$$c(h, k) := (\int_{y \ge \omega} d(y)y^{-2(k+1)}h(y)) + \sum_{n \ge -k} (j_{n+k}/(2n-1)) \cdot \omega^{2n-1}.$$

h₁(z) is an entire function of z on C: $h_1(z) = \sum_{k \ge 0} (z^{2k+1})(-1)^k \cdot c(h, k)$. Therefore $h^\#(z) = h_1(z) + (\pi/2)h(iz)$ yields the analytic extension of h(z) from C − iωL to C.

$$h^{\#}(z) = ( \sum_{k \geq 0} (z^{2k+1})(-1)^k \cdot c(h, k)) + (\pi/2)h(iz).$$

Thus **Lemma 1** holds.

Note that $c(h, k) = \lim_{\omega > 0, \omega \to \infty} \sum_{n \geq 1} (j_{n+k}/(2n-1)) \cdot \omega^{2n-1}$.

The following **Conditional corollary 1** is obtained from **Lemma 1**.

**Review Part I**, **§6**, definitions of j(u) and υ(z).

**Conditional corollary 1** *Assume ('): $\theta < 1$ and for $v \geq 0$, $|P_0(v)| \sim O((log(v))^{\theta})$, as $v \to \infty$.*
*(1) $(z/\pi)\int_{y \geq 0} d(y)(1/(y^2 + z^2))j(y)$ converges absolutely to an analytic function on the half-plane $Re(z) > 0$.*
**Definition of υ(z).** *Set $υ(z) := (1/\pi)j^{\#}(z)$, when $Re(z) > 0$.*
*(2) υ(z) extends to an entire function on C.*

**Proof** Apply **Lemma 1** with h(y) = j(y) .

**Review Part I**, **§5**, **Introduction**, **Definition of e(z)**.

**Conditional corollary 2** *Assume (\*): A is finite and $\lambda(y) = j(y)$, for $y > 0$.*
**(1) $e(-z) = υ(z)$ on the half-plane $Re(z) > 0$:**

$$\sum_{k \geq 1} c(i\gamma e_k)exp(-\gamma_k z) = (z/\pi)\int_{y > 0} d(y)(1/(z^2 + y^2))(-c(0) + P_0(\pi e^{2y}) + P_0(\pi e^{-2y})).$$

**(2) e(s) has an analytic extension from the half-plane $Re(s) < 0$ to the entire complex plane.**

**Proof** Assume A is finite. Say $y \geq 0$. $\lambda(y)$ is bounded. Assume $\lambda(y) = j(y)$. Then (') of **Lemma 1** holds with h(y) = j(y). One obtains (1), (2) of the previous **Conditional corollary 1**. $A < \infty$ together with the **Conditional lemma 6.1** established in **Part I**, **§6**, gives e(-z) =$(1/\pi)\lambda^{\#}(z)$, when $Re(z) > 0$. $\lambda(y) = j(y)$ yields e(-z) = υ(z). **Conditional corollary 1 (2)** now reveals that υ(-s) is the analytic extension of e(s) to all of C.

The proof of the **Main conditional theorem (1) (i), (i')**, was completed in **Part IV**. That now sparks the genesis of the following **Conditional corollaries 3-4**.

Each $\gamma_n$, $\zeta'(\frac{1}{2} + i\gamma_n)$ of the sequence $\gamma_1$, $\zeta'(\frac{1}{2} + i\gamma_1)$; ...$\gamma_n$, $\zeta'(\frac{1}{2} + i\gamma_n)$; ... can successively be expressed in terms of the predecessors $\gamma_k$, $\zeta'(\frac{1}{2} + i\gamma_k)$, with $1 \leq k$ $\leq$ n − 1, and, in the case of $\zeta'(\frac{1}{2} + i\gamma_n)$, also $\gamma_n$.

(See **Part I, §6, Conditional corollary 6.3**)

**Conditional corollary 3** *Assume C^.*
**(1)** *$e(-z) = v(z)$, provided Re(z) > 0. v(z) extends e(-z) to an entire function on C.*
**(2)** *Relations of $\gamma_n$, $\zeta'(\frac{1}{2} + i\gamma_n)$ to their predecessors and the $\zeta(\frac{1}{2} + 4k)$.*

$$\gamma_1 = -lim_{x > 0,\, x \to \infty}\ (1/x)log((-1)^n \cdot\cdot v(x)).$$

$$\zeta'(\frac{1}{2} + i\gamma_1) = 1/(b(i\gamma_1))lim_{Re(z) \to \infty}\ exp(\gamma_1 z)\cdot v(z).$$

$$\gamma_n = -lim_{x > 0,\, x \to \infty}\ (1/x)log((-1)^n\ (v(x) − e(-x, n − 1))).$$

$$\zeta'(\frac{1}{2} + i\gamma_n) = 1/(b(i\gamma_n))lim_{Re(z) \to \infty}\ exp(\gamma_n z)(v(z) − e(-z, n − 1))).$$

## §2 Representation of $p_{i,+}(z)$ via j(y).

**Review** Part I, §6. Definitions of $p_{i,+}(z)$, $\Theta(\theta, z)$. Conditional corollaries 6.4-6.5. Conditional corollary 6.5 therein states without proof the following.

**Conditional corollary 4  Representation of $p_{i,+}(z)$ via j(y).**
(1) *Assume A is finite and for y > 0, $\lambda(y) = j(y)$. Say Im(z) < 0. Then*

$$p_{i,+}(z) = \int_{\theta > 0} d(\theta)j(\theta)\Theta(\theta, z).$$

(2) *Assume C^. Then the previous representation of $p_{i,+}(z)$ holds on the lower half-plane of z with Im(z) < 0.*

**Proof of (1).** Assume A is finite and Re(s) < $\gamma_1$. Then Conditional corollary 6.4, proven in §6 of Part I, gives

$$ip_{i,+}(is) = \int_{y > 0} d(y)e^{sy}(-e(-y)).$$

Also assume $\lambda(y) = j(y)$, for y > 0. Then $|j(y)| \leq A$. Conditional corollary 2 of §1 yields

$$e(-u) = (u/\pi)\int_{\theta \geq 0} d(\theta)(1/(\theta^2 + u^2))j(\theta),$$

when Re(u) > 0. Then

$$ip_{i,+}(is) = \int_{y>0} d(y)e^{sy}(-y/\pi)\int_{\theta \ge 0} d(\theta)(1/(\theta^2 + y^2))j(\theta).$$

Set $x := \mathrm{Re}(s)$. Say $x < 0$. Then

$$\int_{y>0} d(y)\int_{\theta \ge 0} d(\theta)|e^{sy}(-y/\pi)(1/(\theta^2 + y^2))j(\theta)| \le \alpha'\alpha A.$$

Here $\alpha := (y/\pi)\int_{\theta \ge 0} d(\theta)(1/(\theta^2 + y^2)) = \frac{1}{2}$ and $\alpha' := \int_{y>0} d(y)e^{xy} = -1/x$. Apply the interchange $\int_{y>0} d(y)\int_{\theta \ge 0} d(\theta) = \int_{\theta \ge 0} d(\theta)\int_{y>0} d(y)$.

**Proof of (2).** Assume $C^{\wedge}$. Then A is finite. Next apply the Main conditional theorem (1) (i), stated in Part I, §5, (5.4), and proven in §2 of Part IV. $\lambda(y) = j(y)$ results. Thus (2) of Conditional corollary 4 follows from (1) thereof.